# Generalized Reflection Coefficients in Toeplitz-Block-Toeplitz Matrix Case and Fast Inverse 2D levinson Algorithm

April 4, 2004

Rami Kanhouche[1]

## Abstract

A factorization of the inverse of a Hermetian positive definite matrix based on a diagonal by diagonal recurrence formulae permits the inversion of Toeplitz Block Toeplitz matrices using minimized matrix-vector products, with a complexity of $O\left(n_1^3 n_2^2\right)$, where $n_1$, is the block size, and $n_2$ is the block matrix size. A 2D levinson algorithm is introduced that outperform Wittle, Wiggins and Robinson Algorithm.

## Introduction.

Fast inversion algorithms for Toeplitz forms are of important role in our days, these represent also the mathematical solutions to many different problems arising in signal processing theory. In this work we treat the case of linear system with Hermetian Toeplitz Block Toeplitz Matrix, which is, with a direct connection to the 2D Autoregressive model. In the case of Hermetian Toeplitz matrix the famous levinson algorithm [6], first introduced the concept of reflection coefficients, this was further generalized in [1] into the generalized reflection coefficients which are applied to any Strongly Regular Hermetian matrix, in [2] the special case of Hermetian Block Toeplitz was treated and an algorithm was established. In this paper we treat the case of positive definite Hermetian TBT matrix, The algorithm that will be introduced in [2] and the one that will be in this paper, provide a better counterparts in calculations coast for the one presented by Wittle, Wiggins and Robinson [7][8], and its 2D version, respectively. For a better overview of the later algorithms with connection to 2D signal processing the reader is invited to examine [3]. Even that the inverse is factorized, for the benefit of faster application; using Gohberg-Semencul Formula [4][5], would be most appropriate to fully calculate the inverse, and fully taking advantage of the current work.

## Notations.

We consider the Hermetian positive definite (covariance) matrix,

$$R := \left\{ r_{i,j} \right\}_{i,j:=0,1,\ldots n-1} \qquad (1)$$

$r_{i,j} \in C$

$n \in N$

For each matrix *R* we consider the sub matrices family:

[1] PhD Student at Lab. CMLA, Ecole Normale Supérieure de Cachan, 61, avenue du Président Wilson, 94235 CACHAN Cedex, France. phone: +33-1-40112688, mobile: +33-6-62298219, fax: +33-1-47405901, e-mail: rami.kanhouche@cmla.ens-cachan.fr, kanram@wanadoo.fr.

$$R^{k,l} := \left\{ b_{i,j} \right\}_{i,j:=0,\ldots l-k}, b_{i,j} := r_{i+k,j+k} \quad \begin{matrix} k < l, \\ k,l := 0\ldots n-1 \end{matrix}$$

For each $R^{k,l}$ matrix we define the following block diagonal matrix:

$$^{0}R^{k,l} := \begin{bmatrix} 0_{k\times k} & & \\ & R^{k,l} & \\ & & 0_{(n-l)\times(n-l)} \end{bmatrix}$$

which contains $R^{k,l}$ diagonally positioned, with zero filled elsewhere.
The transport square matrix $J$ defined as:

$$J := \begin{bmatrix} 0 & \cdots & 0 & 1 \\ \vdots & \iddots & 1 & 0 \\ 0 & \iddots & \iddots & \vdots \\ 1 & 0 & \cdots & 0 \end{bmatrix}$$

## 1 Reflection Coefficients and Choleskey Factorization of inverse for Hermetian Positive definite matrices.

For any Hermetian p.d matrix the inverse can be factorized into the form [1]

$$R^{-1} = \left(R^{P}\right)^{*} \left(D^{P}\right)^{-1} \left(R^{P}\right)^{T} \quad (2)$$

where $D^{P} = \left(R^{P}\right)^{H} (R)\left(R^{P}\right)$ (3)

is diagonal, and $R^{P}$ is a lower triangular matrix formed from Vectors $\mathrm{p}_{k,l}$, $k,l = 0,1,\ldots n-1$ , $k \le l$ according to

$$R^{P} = \left[\mathrm{p}_{0,n-1}, \mathrm{p}_{1,n-1}, \ldots, \mathrm{p}_{n-1,n-1}\right] \quad (4)$$

These vectors (or polynomials) are mutually orthogonal relatively to the matrix defined cross product:

$$\langle v_1, v_2 \rangle := v_1^{T} R v_2$$

$$\forall\, 0 \le k, 0 \le k', k < l, k' < l, l \le n-1$$

$$\langle \mathrm{p}_{k,l}, \mathrm{p}_{k',l} \rangle = 0 \; if \;\; k \ne k'$$

$$\langle \mathrm{p}_{k,l}, \mathrm{p}_{k',l} \rangle > 0 \; if \;\; k = k'$$

This also do apply for a second family of Polynomials $q_{k,l}$, which also realize the following orthogonally:

$$\forall\, 0 \le k, k < l, k < l', l \le n-1, l' \le n-1$$

$$\langle \mathrm{q}_{k,l'}, \mathrm{q}_{k,l'} \rangle = 0 \; if \;\; l \ne l'$$

$$\langle \mathrm{q}_{k,l'}, \mathrm{q}_{k,l'} \rangle > 0 \; if \;\; l = l'$$

The two families $\mathrm{p}_{k,l}, \mathrm{q}_{k,l}$ are produced from $R$ in a recursive manner. The positive definiteness of the matrix R permits us to establish a recurrence system in which the famous reflection coefficients are used. The Generalized reflection coefficients define the relation between a group of vectors $\mathrm{p}_{k,k+d}, \mathrm{q}_{k,k+d}$ and the next step group of vectors $\mathrm{p}_{k,k+d+1}, \mathrm{q}_{k,k+d+1}$, according to

$$\mathrm{p}_{k,l} = \mathrm{p}_{k,l-1} - a_{k,l}\, \mathrm{q}_{k+1,l} \quad (5)$$

$$\mathrm{q}_{k,l} = \mathrm{q}_{k+1,l} - a'_{k,l}\, \mathrm{p}_{k,l-1} \qquad (6)$$

starting from the canonical basis vectors $\mathrm{p}_{k,k} := \mathrm{e}_k, \mathrm{q}_{k,k} := \mathrm{e}_k \quad k = 0,1,\ldots n-1$ .

The Generalized reflection coefficients are complex numbers, obtained in each step according to:

$$a_{k,l} = \frac{\mathrm{p}_{k,l-1}^T R\, \mathrm{e}_l}{\mathrm{q}_{k+1,l}^T R\, \mathrm{e}_l} \qquad (7)$$

$$a'_{k,l} = \frac{\mathrm{q}_{k+1,l}^T R\, \mathrm{e}_k}{\mathrm{p}_{k,l-1}^T R\, \mathrm{e}_k} \qquad (8)$$

The definition of $v_{k,l} := \mathrm{q}_{k,l}^T R\, \mathrm{e}_l$ , $v'_{k,l} := \mathrm{p}_{k,l}^T R\, \mathrm{e}_k$ install the following recurrence system:

$$v_{k,l} = v_{k+1,l}\left(1 - a_{k,l} a'_{k,l}\right) \in \mathbf{R}^+ \qquad (9)$$

$$v'_{k,l} = v'_{k,l-1}\left(1 - a_{k,l} a'_{k,l}\right) \in \mathbf{R}^+ \qquad (10)$$

This permit us to avoid applying the product in the denominator in (7) and (8) at each step, while for numerators, the following hold:

$$\mathrm{p}_{k,l-1}^T R\, \mathrm{e}_l = \left(\mathrm{q}_{k+1,l}^T R\, \mathrm{e}_k\right)^* \qquad (11)$$

## 2 Generalized reflection coefficients in Hermetian Toeplitz Block Toeplitz Case.

In the Toeplitz-Block-Toeplitz Case each Matrix Block is a Toeplitz matrix of size $n_1$, while the number of Block Matrix is of size $n_2$. We can formalize this structure in the following property.

**Property 1.**
For any $i \le j$

$$r_{i,j} = r_{i \bmod n_1 - j \bmod n_1,\, j \sec n_1 - i \sec n_1} \quad \mathit{if} \quad j \bmod n_1 < i \bmod n_1$$

$$r_{i,j} = r_{0,j-i} \quad \mathit{if} \quad j \bmod n_1 \ge i \bmod n_1$$

*where* $a \sec b := b \operatorname{int}(a/b)$, $\operatorname{int}(x)$ *is the integral part of x.*

*Proof.* From [2], property 1, according to the Block-Toeplitz property we can write,
*for* $i \le j$, $r_{i,j} = r_{i \bmod n_1,\, j - i \sec n_1}$
in the case of the Blocks being Toeplitz matrices, for each element inside the block we differentiate between two cases:

i) $(j - i \sec n_1) \bmod n_1 < i \bmod n_1$: Element exclusively in the lower triangular part of the block.

ii) $(j - i \sec n_1) \bmod n_1 \ge i \bmod n_1$: Element in the upper triangular part of the block.

From this, one can project back on the first line and first column elements of each block in the first block line:
*For* $i \le j$

$$r_{i,j} = r_{i \bmod n_1 - (j - i \sec n_1) \bmod n_1,\, j - i \sec n_1 - (j - i \sec n_1) \bmod n_1} \quad \mathit{if} \quad (j - i \sec n_1) \bmod n_1 < i \bmod n_1$$

$$r_{i,j} = r_{0,(j - i \sec n_1) - i \bmod n_1} \quad \mathit{if} \quad (j - i \sec n_1) \bmod n_1 \ge i \bmod n_1 .$$

While for $i > j$ , $r_{i,j} = \left(r_{j,i}\right)^{*}$.
From which we obtain our proof.□

**Theorem 1.** *In the case of Toeplitz block Toeplitz Matrix the generalized reflection coefficients admit the following:*

$$a_{k,l} = a'^{*}_{k',l'}$$

$$a'_{k,l} = a^{*}_{k',l'}$$

$$\mathrm{p}_{k,l} = U^{k-k'}\left(\tilde{\mathrm{q}}_{k',l'}\right)^{*}$$

$$\mathrm{q}_{k,l} = U^{k-k'}\left(\tilde{\mathrm{p}}_{k',l'}\right)^{*}$$

$$v_{k,l} = v'^{*}_{k',l'}$$

$$v'_{k,l} = v^{*}_{k',l'}$$

*where* , $\left(k',l'\right) := \left(k \sec n_1 + n_1 - 1 - l \bmod n_1, l \sec n_1 + n_1 - 1 - k \bmod n_1\right)$, *and for any*

$$v_{k,l} := \left[x_i\right]_{i:=0,..w}, \quad k,l \le w \ , \ \tilde{v}_{k,l} := \left[\tilde{x}_i\right]_{i:=0,..w}, \begin{cases} \tilde{x}_i := x_{l-i+k} & if \ k \le i \le l \\ \tilde{x}_i := x_i \ if \ i < k, or \ i > l \end{cases}$$

$U$ is the $\left(n_1 n_2\right)\times\left(n_1 n_2\right)$ up-down shift matrix defined as

$$U := \begin{bmatrix} 0 & 0 & \cdots & 0 & 0 \\ 1 & 0 & \ddots & & 0 \\ 0 & \ddots & \ddots & 0 & \vdots \\ \vdots & \ddots & 1 & 0 & 0 \\ 0 & \cdots & 0 & 1 & 0 \end{bmatrix}, U^{-1} := \begin{bmatrix} 0 & 1 & 0 & \cdots & 0 \\ 0 & 0 & \ddots & \ddots & \vdots \\ 0 & \ddots & \ddots & 1 & 0 \\ \vdots & \ddots & 0 & 0 & 1 \\ 0 & \cdots & 0 & 0 & 0 \end{bmatrix} \text{ and } U^{0} := \mathbf{1} \text{ is the identity}$$

matrix.

The following properties will help us establishing our proof.
We will note $\mathbf{p}_{k,l}, \mathbf{q}_{k,l}, a_{k,l}, a'_{k,l}, v_{k,l}, v'_{k,l}$ , $0 \le k < l \le n-1$, corresponding to a Hermetian matrix $M$ of size $\left(n \times n\right)$, respectively as $\mathrm{p}^{M}_{k,l}, \mathrm{q}^{M}_{k,l}, a^{M}_{k,l}, a'^{M}_{k,l}, v^{M}_{k,l}, v'^{M}_{k,l}$ .

**Lemma 1.**

$\forall h < s, s : s \bmod r = 0$

i) $\left(s-h\right) \bmod r = r - 1 - \left(h-1\right) \bmod r$

ii) $\left(s-h\right) \sec r = s - r - \left(h-1\right) \sec r$

*Proof.* For any $\left(s-h\right) \bmod r = k$ , it is clear that we can write $h$ in the form:
$h = tr + r - k, \quad t \in N$ . Also it is clear that

$\forall x, r, w : x \bmod r = w$

$w > n \Rightarrow \left(x-n\right) \bmod r + n = w$

Since for $k = 1,2,\ldots r-1$:

$h \bmod r = r - k$ , this also realize $\left(h-1\right) \bmod r + 1 = r - k$ , which can be reorganized to yield $k = r - \left(h-1\right) \bmod r - 1$, the equality also hold for the case of $k = 0$. Finally relation (ii) is direct replacement of (i) result in the definition of *sec* according to $x \sec r = x - x \bmod r$ .□

**Property 2.** *In the case of Hermetian Toeplitz-block-Toeplitz matrix the sub matrices admit the following:*

$R_{k,l} = J\left(R_{k\otimes l}\right)^T J$ *, where* $\left(k \otimes l\right) := \left(k', l'\right)$

*Proof.* We start by noting the elements of each sub matrix $M$ as $r_{i,j}^{M}$ , from which we got $r_{i,j}^{R_{k,l}} = r_{i+k,j+k}$ , where $i, j = 0,1,\ldots,l-k$ . (1)

this can be put into the form $r_{i,j}^{R_{k,l}} = r_{f_1(k,i),f_1(k,j)}$, where $f_1\left(k,t\right) = k + t$ .

Next we notice that for any two matrices $M_1, M_2$ of size *n*:

$M_1 = JM_2 \Leftrightarrow r_{i,j}^{M_1} = r_{n-1-i,j}^{M_2}$

$M_1 = M_2 J \Leftrightarrow r_{i,j}^{M_1} = r_{i,n-1-j}^{M_2}$

from which we get :

$M_1 = M_2^T J \Leftrightarrow r_{i,j}^{M_1} = r_{j,n-1-i}^{M_2}$

and $M_1 = JM_2^T J \Leftrightarrow r_{i,j}^{M_1} = r_{n-1-i,j}^{M_2^T J} \Leftrightarrow r_{i,j}^{M_1} = r_{n-1-i,n-1-j}^{M_2^T} \Leftrightarrow r_{i,j}^{M_1} = r_{n-1-j,n-1-i}^{M_2}$ .

and eventually $M_1 = JM_2^T J \Leftrightarrow r_{i,j}^{M_1} = r_{n-1-j,n-1-i}^{M_2}$ .

Or $r_{i,j}^{J\left(R_{k\otimes l}\right)^T J} = r_{(l-k)-j,(l-k)-i}^{R_{k\otimes l}}$ , from the definition of $k \otimes l$ , and (1),

$r_{i,j}^{R_{k\otimes l}} = r_{(k \sec n_1 + n_1 - 1 - l \bmod n_1)+i,(k \sec n_1 + n_1 - 1 - l \bmod n_1)+j}$ which give

$r_{i,j}^{J\left(R_{k\otimes l}\right)^T J} = r_{(k \sec n_1 + n_1 - 1 - l \bmod n_1)+(l-k)-j,(k \sec n_1 + n_1 - 1 - l \bmod n_1)+(l-k)-i}$

and simplified into

$r_{i,j}^{J\left(R_{k\otimes l}\right)^T J} = r_{f_2(k,l,j),f_2(k,l,i)}$ , where $f_2\left(k,l,t\right) = l \sec n_1 - k \bmod n_1 + n_1 - 1 - t$

According to *Lemma 1*, by replacing with $r = n_1, s = l \sec n_1 + n_1$ , and $h = k \bmod n_1 + 1 + t$ ,

$\left(f_1\left(k,t\right)\right) \bmod n_1 = n_1 - 1 - \left(f_2\left(k,l,t\right)\right) \bmod n_1$

$\left(f_1\left(k,t\right)\right) \sec n_1 = l \sec n_1 - \left(f_2\left(k,l,t\right)\right) \sec n_1$

this will oblige the following

$\left(f_1\left(k,i\right)\right) \bmod n_1 > \left(f_1\left(k,j\right)\right) \bmod n_1 \Leftrightarrow \left(f_2\left(k,l,i\right)\right) \bmod n_1 < \left(f_2\left(k,l,j\right)\right) \bmod n_1$

$\left(f_1\left(k,i\right)\right) \bmod n_1 - \left(f_1\left(k,j\right)\right) \bmod n_1 = \left(f_2\left(k,l,j\right)\right) \bmod n_1 - \left(f_2\left(k,l,i\right)\right) \bmod n_1$

$\left(f_1\left(k,j\right)\right) \sec n_1 - \left(f_1\left(k,i\right)\right) \sec n_1 = \left(f_2\left(k,l,i\right)\right) \sec n_1 - \left(f_2\left(k,l,j\right)\right) \sec n_1$

Since $f_2\left(k,l,i\right) - f_2\left(k,l,j\right) = f_1\left(k,j\right) - f_1\left(k,i\right)$ , we did cover both conditions of *property 1*, in the case of $f_1\left(k,i\right) \le f_1\left(k,j\right)$, while for the case $f_1\left(k,i\right) > f_1\left(k,j\right)$ the equality of matrices $R_{k,l}, J\left(R_{k\otimes l}\right)^T J$ follows from each being Hermetian, with which the proof is completed.□

**Property 3.** *For any Hermetian matrix M of size* $\left(n \times n\right)$ *,*

$$\left(\mathrm{p}_{k,l}^{JM^TJ},\mathrm{q}_{k,l}^{JM^TJ},a_{k,l}^{JM^TJ},a_{k,l}^{JM^TJ},v_{k,l}^{JM^TJ},v'^{JM^TJ}_{k,l}\right)=$$

$$\left(\left(\tilde{\mathrm{q}}_{n-1-l,n-1-k}^{M}\right)^*,\left(\tilde{\mathrm{p}}_{n-1-l,n-1-k}^{M}\right)^*,\left(a'^{M}_{n-1-l,n-1-k}\right)^*,\left(a^{M}_{n-1-l,n-1-k}\right)^*,\left(v'^{M}_{n-1-l,n-1-k}\right)^*,\left(v^{M}_{n-1-l,n-1-k}\right)^*\right)$$

*Proof.* To proof this we apply the main algorithm (5-10), on the matrix $JM^TJ$.
We notice that, at each step we got-after eliminating *J*:

$$a_{k,l}^{JM^TJ}=\frac{\left(\tilde{\mathrm{p}}_{k,l-1}^{JM^TJ}\right)^T M^T \mathbf{e}_{n-1-l}}{\left(\tilde{\mathbf{q}}_{k+1,l}^{JM^TJ}\right)^T M^T \mathbf{e}_{n-1-l}} \qquad a'^{JM^TJ}_{k,l}=\frac{\left(\tilde{\mathrm{q}}_{k+1,l}^{JM^TJ}\right)^T M^T \mathbf{e}_{n-1-k}}{\left(\tilde{\mathrm{p}}_{k,l-1}^{JM^TJ}\right)^T M^T \mathbf{e}_{n-1-k}}$$

taking the conjugate on both sides of each relation, since $M^T=M^*$, we get

$$\left(a_{k,l}^{JM^TJ}\right)^*=\frac{\left(\tilde{\mathrm{p}}_{k,l-1}^{JM^TJ}\right)^{*T} M\ \mathrm{e}_{n-1-l}}{\left(\tilde{\mathrm{q}}_{k+1,l}^{JM^TJ}\right)^{*T} M\ \mathrm{e}_{n-1-l}} \qquad \left(a'^{JM^TJ}_{k,l}\right)^*=\frac{\left(\tilde{\mathrm{q}}_{k+1,l}^{JM^TJ}\right)^{T*} M\ \mathbf{e}_{n-1-k}}{\left(\tilde{\mathrm{p}}_{k,l-1}^{JM^TJ}\right)^{T*} M\ \mathbf{e}_{n-1-k}}$$

by defining $s:=n-1-k \quad h:=n-1-l$, We can write the above relations in the form:

$$\left(a_{n-1-s,n-1-h}^{JM^TJ}\right)^*=\frac{\left(\tilde{\mathrm{p}}_{n-1-s,n-1-(h+1)}^{JM^TJ}\right)^{*T} M\ \mathrm{e}_h}{\left(\tilde{\mathrm{q}}_{n-1-(s-1),n-1-h}^{JM^TJ}\right)^{*T} M\ \mathrm{e}_h} \qquad \left(a'^{JM^TJ}_{n-1-s,n-1-h}\right)^*=\frac{\left(\tilde{\mathrm{q}}_{n-1-(s-1),n-1-h}^{JM^TJ}\right)^{T*} M\ \mathrm{e}_s}{\left(\tilde{\mathrm{p}}_{n-1-s,n-1-(h+1)}^{JM^TJ}\right)^{T*} M\ \mathrm{e}_s}$$

Next we take advantage of the definitions:

$$\vec{\mathrm{p}}_{k,l}:=\left(\tilde{\mathrm{p}}_{n-1-l,n-1-k}^{JM^TJ}\right)^* \qquad \vec{\mathrm{q}}_{k,l}:=\left(\tilde{\mathrm{q}}_{n-1-l,n-1-k}^{JM^TJ}\right)^*$$

$$\vec{a}_{k,l}:=\left(a_{n-1-l,n-1-k}^{JM^TJ}\right)^* \qquad \vec{a}'_{k,l}:=\left(a'^{JM^TJ}_{n-1-l,n-1-k}\right)^*$$

our definitions hold for all values $k\le l,\quad k,l=0,\ldots n-1$. So we can write in the form:

$$\vec{a}_{h,s}=\frac{\left(\vec{\mathrm{p}}_{h+1,s}\right)^T M\ \mathbf{e}_h}{\left(\vec{\mathrm{q}}_{h,s-1}\right)^T M\ \mathbf{e}_h} \qquad \vec{a}'_{h,s}=\frac{\left(\vec{\mathrm{q}}_{h,s-1}\right)^T M\ \mathbf{e}_s}{\left(\vec{\mathrm{p}}_{h+1,s}\right)^T M\ \mathbf{e}_s} \quad \text{(p3.1)}$$

also by taking the conjugate on both sides of (5-6),and multiplying by *J* we obtain :

$$\vec{\mathrm{p}}_{h,s}=\vec{\mathrm{p}}_{h+1,s}-\left(a_{n-1-s,n-1-h}^{JM^TJ}\right)^*\vec{\mathrm{q}}_{h,s-1} \qquad \vec{\mathrm{q}}_{h,s}=\vec{\mathrm{q}}_{h,s-1}-\left(a'^{JM^TJ}_{n-1-s,n-1-h}\right)^*\vec{\mathrm{p}}_{h+1,s}$$

which yield

$$\vec{\mathrm{p}}_{h,s}=\vec{\mathrm{p}}_{h+1,s}-\vec{a}_{h,s}\,\vec{\mathrm{q}}_{h,s-1} \qquad \vec{\mathrm{q}}_{h,s}=\vec{\mathrm{q}}_{h,s-1}-\vec{a}'_{h,s}\,\vec{\mathrm{p}}_{h+1,s} \quad \text{(p3.2)}$$

comparing the new relations (p3.1,p3.2) with the basic algorithm we conclude that $\vec{\mathrm{p}}_{k,l}=\mathrm{q}_{k,l}^M\ \vec{\mathrm{q}}_{k,l}=\mathrm{p}_{k,l}^M$ and $\vec{a}_{k,l}=a'^M_{k,l}\quad \vec{a}'_{k,l}=a^M_{k,l}$, from which our proof is completed.□

**Property 4.**

$$\left(\mathrm{p}_{s,d}^{{}^0R^{k,l}},\mathrm{q}_{s,d}^{{}^0R^{k,l}},a_{s,d}^{{}^0R^{k,l}},a'^{{}^0R^{k,l}}_{s,d},v_{s,d}^{{}^0R^{k,l}},v'^{{}^0R^{k,l}}_{s,d}\right)=\left(\mathrm{p}_{s,d}^R,\mathrm{q}_{s,d}^R,a_{s,d}^R,a'^R_{s,d},v_{s,d}^R,v'^R_{s,d}\right)$$

With condition that $k\le s<d\le l$.
Where $\mathrm{p}_{k,l},\mathrm{q}_{k,l},a_{k,l},a'_{k,l},v_{k,l},v'_{k,l}$ that corresponds to a given covariance matrix *M*, are noted respectively as $\mathrm{p}_{k,l}^M,\mathrm{q}_{k,l}^M,a_{k,l}^M,a'^M_{k,l},v_{k,l}^M,v'^M_{k,l}$.
The property is so obvious and need no proof.

**Proof Of Theorem 1.** We start by noticing that $\forall 0 \le k < l \le n_1 n_2 - 1 \; \exists s,d : k \le s < d \le l$ , so that we can write

$$\left(\mathrm{p}_{s,d}^{{}^0R^{k,l}}, \mathrm{q}_{s,d}^{{}^0R^{k,l}}, a_{s,d}^{{}^0R^{k,l}}, a_{s,d}'^{{}^0R^{k,l}}, v_{s,d}^{{}^0R^{k,l}}, v_{s,d}'^{{}^0R^{k,l}}\right) =$$
$$\left(U_{l-k}^{k}\, \mathrm{p}_{s-k,d-k}^{R^{k,l}}, U_{l-k}^{k}\, \mathrm{q}_{s-k,d-k}^{R^{k,l}}, a_{s-k,d-k}^{R^{k,l}}, a_{s-k,d-k}'^{R^{k,l}}, v_{s-k,d-k}^{R^{k,l}}, v_{s-k,d-k}'^{R^{k,l}}\right)$$

where $U_x$ is a Matrix of size $(n_1 n_2) \times (x+1)$, defined as:

$$U_x := \begin{bmatrix} 0 & 0 & \cdots & 0 \\ 1 & 0 & \ddots & \vdots \\ 0 & 1 & \ddots & 0 \\ 0 & 0 & \ddots & 0 \\ \vdots & \vdots & \ddots & 1 \\ 0 & 0 & \cdots & 0 \end{bmatrix}, \quad U_x^0 := \begin{bmatrix} \mathbf{1}_{(x+1)\times(x+1)} \\ \mathbf{0} \end{bmatrix}$$

coupling this with *property 4*, we obtain that

$$\left(\mathrm{p}_{s,d}^{R}, \mathrm{q}_{s,d}^{R}, a_{s,d}^{R}, a_{s,d}'^{R}, v_{s,d}^{R}, v_{s,d}'^{R}\right) =$$
$$\left(U_{l-k}^{k}\, \mathrm{p}_{s-k,d-k}^{R^{k,l}}, U_{l-k}^{k}\, \mathrm{q}_{s-k,d-k}^{R^{k,l}}, a_{s-k,d-k}^{R^{k,l}}, a_{s-k,d-k}'^{R^{k,l}}, v_{s-k,d-k}^{R^{k,l}}, v_{s-k,d-k}'^{R^{s,d}}\right)$$

we select $d = l, s = k$

$$\left(\mathbf{p}_{k,l}^{R}, \mathbf{q}_{k,l}^{R}, a_{k,l}^{R}, a_{k,l}'^{R}, v_{k,l}^{R}, v_{k,l}'^{R}\right) =$$
$$\left(U_{l-k}^{k}\, \mathrm{p}_{0,l-k}^{R^{k,l}}, U_{l-k}^{k}\, \mathrm{q}_{0,l-k}^{R^{k,l}}, a_{0,l-k}^{R^{k,l}}, a_{0,l-k}'^{R^{k,l}}, v_{0,l-k}^{R^{k,l}}, v_{0,l-k}'^{R^{k,l}}\right) \quad \text{(t1.1)}$$

By the same way we obtain for $(k', l') = k \otimes l$

$$\left(\mathrm{p}_{k\otimes l}^{R}, \mathrm{q}_{k\otimes l}^{R}, a_{k\otimes l}^{R}, a_{k\otimes l}'^{R}, v_{k\otimes l}^{R}, v_{k\otimes l}'^{R}\right) =$$
$$\left(U_{l'-k'}^{k'}\, \mathrm{p}_{0,l'-k'}^{R^{k\otimes l}}, U_{l'-k'}^{k'}\, \mathrm{q}_{0,l'-k'}^{R^{k\otimes l}}, a_{0,l'-k'}^{R^{k\otimes l}}, a_{0,l'-k'}'^{R^{k\otimes l}}, v_{0,l'-k'}^{R^{k\otimes l}}, v_{0,l'-k'}'^{R^{k\otimes l}}\right)$$

By combining the application of both *property 2* and *property 3* we obtain

$$\left(\mathrm{p}_{k\otimes l}^{R}, \mathrm{q}_{k\otimes l}^{R}, a_{k\otimes l}^{R}, a_{k\otimes l}'^{R}, v_{k\otimes l}^{R}, v_{k\otimes l}'^{R}\right) =$$
$$\left(U_{l-k}^{k'}\left(\tilde{\mathrm{q}}_{0,l-k}^{R^{k,l}}\right)^{*}, U_{l-k}^{k'}\left(\tilde{\mathrm{p}}_{0,l-k}^{R^{k,l}}\right)^{*}, \left(a_{0,l-k}'^{R^{k,l}}\right)^{*}, \left(a_{0,l-k}^{R^{k,l}}\right)^{*}, \left(v_{0,l-k}'^{R^{k,l}}\right)^{*}, \left(v_{0,l-k}^{R^{k,l}}\right)^{*}\right) \quad \text{(t1.2)}$$

By comparing first sides of each of (t1.1), and (t1.2), with the second sides, we find that

$$a_{k,l} = a_{k',l'}^{*}$$
$$a_{k,l}' = a_{k',l'}'^{*}$$
$$v_{k,l} = v_{k',l'}'^{*}$$
$$v_{k,l}' = v_{k',l'}^{*}$$

To complete our proof it remains to notice that

$$\mathrm{p}_{k,l}^{R} = U_{l-k}^{k}\, \mathrm{p}_{0,l-k}^{R^{k,l}} \quad \text{(t1.3)}$$
$$\mathrm{q}_{k\otimes l}^{R} = U_{l-k}^{k'}\left(\tilde{\mathrm{p}}_{0,l-k}^{R^{k,l}}\right)^{*}$$

by multiplying both sides of the last relation with $\left(U_{l-k}^{k'}\right)^T$, we obtain,

$\left(U_{l-k}^{k'}\right)^T \mathrm{q}_{k\otimes l}^{R} = \left(\tilde{\mathrm{p}}_{0,l-k}^{R^{k,l}}\right)^*$ (t1.4) , since

$\left(U_{l-k}^{k'}\right)^T \left(U_{l-k}^{k'}\right) = \mathbf{1}$. To prove the later relation one should only pay attention that by definition $k' < l' \le n_1 n_2 - 1 \Rightarrow l' - k' \le n_1 n_2 - 1 - k'$, from which we have $k' \le n_1 n_2 - 1 - (l' - k')$, also that the equality $l - k = l' - k'$ is always realized.

From (t1.4) we can further write $\left(U_{l-k}^{k'}\right)^T \left(\tilde{\mathrm{q}}_{k\otimes l}^{R}\right)^* = \mathrm{p}_{0,l-k}^{R^{k,l}}$, by replacing $\mathrm{p}_{0,l-k}^{R^{k,l}}$ from the last relation in (t1.3) we conclude that: $\mathrm{p}_{k,l}^{R} = U_{l-k}^{k} \left(U_{l-k}^{k'}\right)^T \left(\tilde{\mathrm{q}}_{k\otimes l}^{R}\right)^*$. Since, for any vector $v_{k',l'}$ defined as

$v_{b,g} := [x_i]_{i:=0,\ldots n-1} : 0 \le b < g \le n-1, x_i = 0 \; if \;\; i < b \; or \; i > g$

we got always $\forall k \ge 0, U_{g-b}^{k} \left(U_{g-b}^{b}\right)^T v_{b,g} = U^{k-b} v_{b,g}$

we get $\mathrm{p}_{k,l}^{R} = U^{k-k'} \left(\tilde{\mathrm{q}}_{k\otimes l}^{R}\right)^*$, by the same way we obtain that

$\mathrm{q}_{k,l}^{R} = U^{k-k'} \left(\tilde{\mathrm{p}}_{k\otimes l}^{R}\right)^*$. By this the proof is completed.□

The next theorem will help us obtaining our inverse factorization with a minimum of calculus effort, it follows the main steps of the one presented in [2], with the exception of the algorithm being working on the first skew symmetric part of each Toeplitz block.

**Theorem 2.** *In Hermetian Toeplitz Block Toeplitz case, the application of Recurrence relations (5)-(10) resume in the following Algorithm, with complexity* $O\left(n_1^3 n_2^2\right)$*:*
*First we rewrite Equations (5)-(10) as the subroutine,*

*Subroutine 1*

$$a_{k,l} = \frac{\hat{\mathrm{p}}^T R \, \mathrm{e}_l}{\hat{\mathrm{v}}} \qquad a'_{k,l} = \frac{\hat{\mathrm{q}}^T R \, \mathrm{e}_k}{\hat{\mathrm{v}}'} \qquad \text{(t2.1)}$$

$$\mathrm{p}_{k,l} = \hat{\mathrm{p}} - a_{k,l}\,\hat{\mathrm{q}} \qquad \mathrm{q}_{k,l} = \hat{\mathrm{q}} - a'_{k,l}\,\hat{\mathrm{p}} \qquad \text{(t2.1)}$$

$$v_{k,l} = \hat{\mathrm{v}}\left(1 - a_{k,l} a'_{k,l}\right) \qquad v'_{k,l} = \hat{\mathrm{v}}'\left(1 - a_{k,l} a'_{k,l}\right) \qquad \text{(t2.3)}$$

*and proceed as the following,*

*Initialization:*
*For* $k = 0 \; to \; n_1 - 1 : \{$

$\mathrm{p}_{k,k} = \mathrm{q}_{k,k} = e_k$, $v_{k,k} = v'_{k,k} = r_{0,0}$

$\}$

*Main routine:*

*For* $d_2 = 0 \, to \; n_2 - 1 \, do$ :

$\{$

*If* $d_2$ *not equal to zero:*

```
{
    For d_1 = n_1 - 1 to 0 do :  (Lower triangle loop)
    {
        For u = 0 to n_1 - d_1 - 1 do :
        {
        k = u + d_1    l = d_2 n_1 + u
        k' = k sec n_1 + n_1 - 1 - l mod n_1
                if k ≤ k' {
                (k^0, l^-) = (k, l - 1)
                        if u equal to 0:{
                        (k^0', l^-') = (k^0 ⊗ l^-)
                        p̂ = U^{n_1-1} (q̃_{k^0', l^-'})^* , v̂' = v_{k^0', l^-'}  t2.4
                        }
                        else {
                        p̂ = p_{k^0, l^-}   , v̂' = v'_{k^0, l^-}
                        }
                        if k equal to k' :{
                        q̂ = U (p̂~)^* ,  v̂ = v̂'  t2.5
                        }
                        else {
                        (k^+, l^0) = (k + 1, l)
                        q̂ = q_{k^+, l^0}   , v̂ = v_{k^+, l^0}
                        }
                        Apply Subroutine 1
                }
        }
    }
}
For d_1 = 1 to n_1 - 1 do : (Exclusively Upper triangle loop)
{
    For u = 0 to n_1 - d_1 - 1 do :
    {
        k = u    l = d_2 n_1 + u + d_1
        k' = k sec n_1 + n_1 - 1 - l mod n_1
                if k ≤ k' {
                (k^0, l^-) = (k, l - 1)
                p̂ = p_{k^0, l^-}   , v̂' = v'_{k^0, l^-}
                        if k equal to k' :{
                        q̂ = U (p̂~)^* , v̂ = v̂'  t2.6
                        }
                        else {
```

$(k^+, l^0) = (k+1, l)$

$\hat{\text{q}} = \text{q}_{k^+,l^0} \quad, \hat{\text{v}} = v_{k^+,l^0}$

}

*Apply Subroutine 1*

}

}

}

}

*Proof.* For any element $(i,j) \in \text{T}_d$ ,

$\text{T}_d := \{(k,l) : k < l, (k,l) \in [0, n_1 d - 1] \times [0, n_1 d - 1]\}, d = 1, 2, \ldots n_2$

$(i \bmod d, j - i \sec d) \in \text{P}_d$

$\text{P}_d := \{(k,l) : k < l, (k,l) \in [0, n_1 - 1] \times [0, n_1 d - 1]\}$

From Theorem 1, by obtaining $\text{p}_{k,l}, \text{q}_{k,l}, a_{k,l}, a'_{k,l}, v_{k,l}, v'_{k,l}$ for all subscripts $(k,l) \in \text{P}_d$ we obtain directly these values for $(k,l) \in \text{T}_d$ .

The algorithm proceed in calculating values for subscripts contained in the following groups, with respective order:

$\Lambda_1, \Lambda_2, \ldots, \Lambda_{n_2}$ , where $\Lambda_d := \text{P}_d \setminus \text{P}_{d-1}$, and $\text{P}_0 := \varnothing$

In *Lower triangle loop* we obtain values for

$\forall (k,l) \in \Lambda_d^-$ , where $\Lambda_d^- := \{(k,l) \in \Lambda_d : k \geq l \bmod n_1\}$

To demonstrate this, first we notice that for

$\Lambda_d^{-,d_1} := \{(k,l) : (k,l) \in \Lambda_d^-, k - l \bmod n_1 = d_1\}$

$$\Lambda_d^- = \bigcup_{d_1=0}^{n_1-1} \Lambda_d^{-,d_1}$$

we got, for $d_1 \in [0, n_1 - 1]$:

$\Lambda_d^{-,d_1} = \{(k,l) : l \bmod n_1 \in [0, n_1 - 1 - d_1], k = l \bmod n_1 + d_1,\}$

By defining $u := l \bmod n_2$ we conclude that we did consider all $(k,l) \in \Lambda_d^-$ .

Using the same logic we conclude that in *Exclusively Upper Triangle Loop* we obtain values for

$\forall (k,l) \in \Lambda_d^+$, where $\Lambda_d^+ := \{(k,l) \in \Lambda_d : k < l \bmod n_1\}$. It is clear that $\Lambda_d = \Lambda_d^- \cup \Lambda_d^+$

In the *EUTL* it is easy to verify that:

$\forall (k,l) \in \Lambda_d^+, \{(k, l-1), (k+1, l)\} \subset \Lambda_d$ ,

while in the *Lower Triangle Loop* this is not always true, more precisely for $k = n_1 - 1$, or/and $l \bmod n_1 = 0$ . Each of $\Lambda_d^+, \Lambda_d^-$ is divided into two skew symmetric parts relative to the antidiagonal, this division come from the fact that $(k,l)$, and $(k \otimes l)$ mark two symmetric points relative to the antidiagonal of each Block, the condition $k \leq k'$ let us treat only the first skew symmetric part, since by *theorem 1* the other parts values need not to be calculated, a special case arise when $k = k'$ in which needed values from the second parts are reintroduced according to *theorem 1* , t2.5,

t2.6. In the case when $(k,l)\in\Lambda_d^-$, when $l \bmod n_1 = 0$, $(k,l-1)$, is in the second skew symmetric part of $\Lambda_{d-1}$ for which again *theorem 1* let us obtain its values by the means of $(k\otimes l-1)$, t2.4.

For complexity estimation we proceed as the following:

Each Entry in the routine will require operations of $O(l-k)$, by noting $c_1$ as the step constant, $opc$ as the total number of operation, the calculus cost take the form of:

$$opc \approx \frac{1}{2}\left[\sum_{d_2=0}^{n_2-1}\sum_{d_1=1}^{n_1-1}\sum_{u=0}^{n_1-d_1-1} c_1\left(n_1 d_2 + d_1\right) + \sum_{d_2=1}^{n_2-1}\sum_{d_1=0}^{n_1-1}\sum_{u=0}^{n_1-d_1-1} c_1\left(n_1 d_2 - d_1\right)\right]$$

which can be expanded into:

$$opc \approx \frac{1}{2}\begin{bmatrix} (n_1-1)c_1\frac{1}{2}(n_1-1)n_1 - c_1\frac{1}{6}(n_1-1)(n_1)(2n_1-1)+ \\ 2c_1 n_1\frac{1}{2}(n_2-1)n_2(n_1-1)^2 - 2c_1 n_1\frac{1}{2}(n_2-1)n_2\frac{1}{2}(n_1-1)n_1 + \\ (n_1-1)c_1 n_1\frac{1}{2}(n_2-1)n_2 \end{bmatrix} \qquad (12)$$

From which we conclude the complexity order of $n_1^3 n_2^2$.

## A Comparison with the WWR Algorithm, In TBT Case.

WWR Algorithm in the TBT case gives a solution to the linear problem:

$$\mathbf{A}_{n_1,n_2}.R_{n_1,n_2} = \mathbf{r}_{n_1,n_2} \qquad (13)$$

$\mathbf{A}_{n_1,n_2} \equiv \left[A_1^{n_2-1}, A_2^{n_2-1}, \ldots, A_{n_2-1}^{n_2-1}\right]$, where $A_k^l$, $k\le l$ is the well know $k$ matrix polynomial coefficient of order $l$.

And $\mathbf{r}_{n_1,n_2} \equiv \left[R_1, R_2, \cdots, R_{n_2-1}\right]$, $R_i \equiv R^{in_1,(i+1)n_1-1}$

$$R_{n_1,n_2} \equiv \begin{bmatrix} R_0 & R_1 & \cdots & R_{n_2-2} \\ R_{-1} & R_0 & \ddots & \vdots \\ \vdots & \ddots & \ddots & R_1 \\ R_{2-n_2} & \cdots & R_{-1} & R_0 \end{bmatrix}$$

The Coefficients Matrix are calculated in a recursive manner, starting from $n>0$, up until $n = n_2 - 1$, according to the order [7][8][9]:

$$\Delta_n \equiv R_n + \sum_{l=1}^{n-1} A_l^{n-1} R_{n-l}$$

$$A_n^n \equiv -\Delta_n\left(P_{n-1}\right)^{-1}$$

$$P_n \equiv P_{n-1} + A_n^n \Delta_n^H$$

With,

$$\mathbf{A}_{n_1,n+1} = \left[\mathbf{A}_{n_1,n-1}, \mathbf{0}\right] + A_n^n\left[JA_{n-1}^{n-1*}J, JA_{n-2}^{n-1*}J, \cdots JA_1^{n-1*}J, \mathbf{1}\right]$$

Starting From

$$P_0 = R_0$$

The total calculus cost of WWR in the TBT case can be written as:

$$opcwwr = \sum_{n=1}^{n_2-1}\left[3(n_1)^3\right] + \sum_{n=2}^{n_2-1}\left[2\sum_{i=1}^{n-1}(n_1)^3\right]$$

$$= (n_2-1)3(n_1)^3 + 2(n_1)^3\sum_{n=2}^{n_2-1}\left[(n-1)\right]$$

$$= (n_2-1)3(n_1)^3 + 2(n_1)^3\frac{\left(4+(n_2-3)\right)(n_2-2)}{2}$$

$$opcwwr = (n_2-1)3(n_1)^3 + (n_1)^3(n_2+1)(n_2-2) \qquad (14)$$

to be able to draw a precise comparison we need to more precise the expression in (12), for that we will set the constant factor $c_1 = 3$, which take into account the Matrix and Vector multiplication operations in (t2.1)(t2.2), while for the divisions in (t2.1), (t2.2), and multiplications in (t2.3) we will add a total value of $\frac{5}{2}(n_1)^2 n_2$ , from which we can write the total number of operation according to *theorem 2* as :

$$opc = (n_1-1)\frac{3}{4}(n_1-1)n_1 - \frac{3}{12}(n_1-1)(n_1)(2n_1-1) +$$
$$n_1\frac{3}{2}(n_2-1)n_2(n_1-1)^2 - n_1(n_2-1)n_2\frac{3}{4}(n_1-1)n_1 + \qquad (15)$$
$$(n_1-1)n_1\frac{3}{4}(n_2-1)n_2 + \frac{5}{2}(n_1)^2 n_2$$

*Graph 1* shows the plot for both relation (14), and (15), for $n_1 = n_2$; WWR Algorithm in dashed line, and *Theorem 2* algorithm in solid line. *Graph 2* plots the ratio for the same values of Graph1.

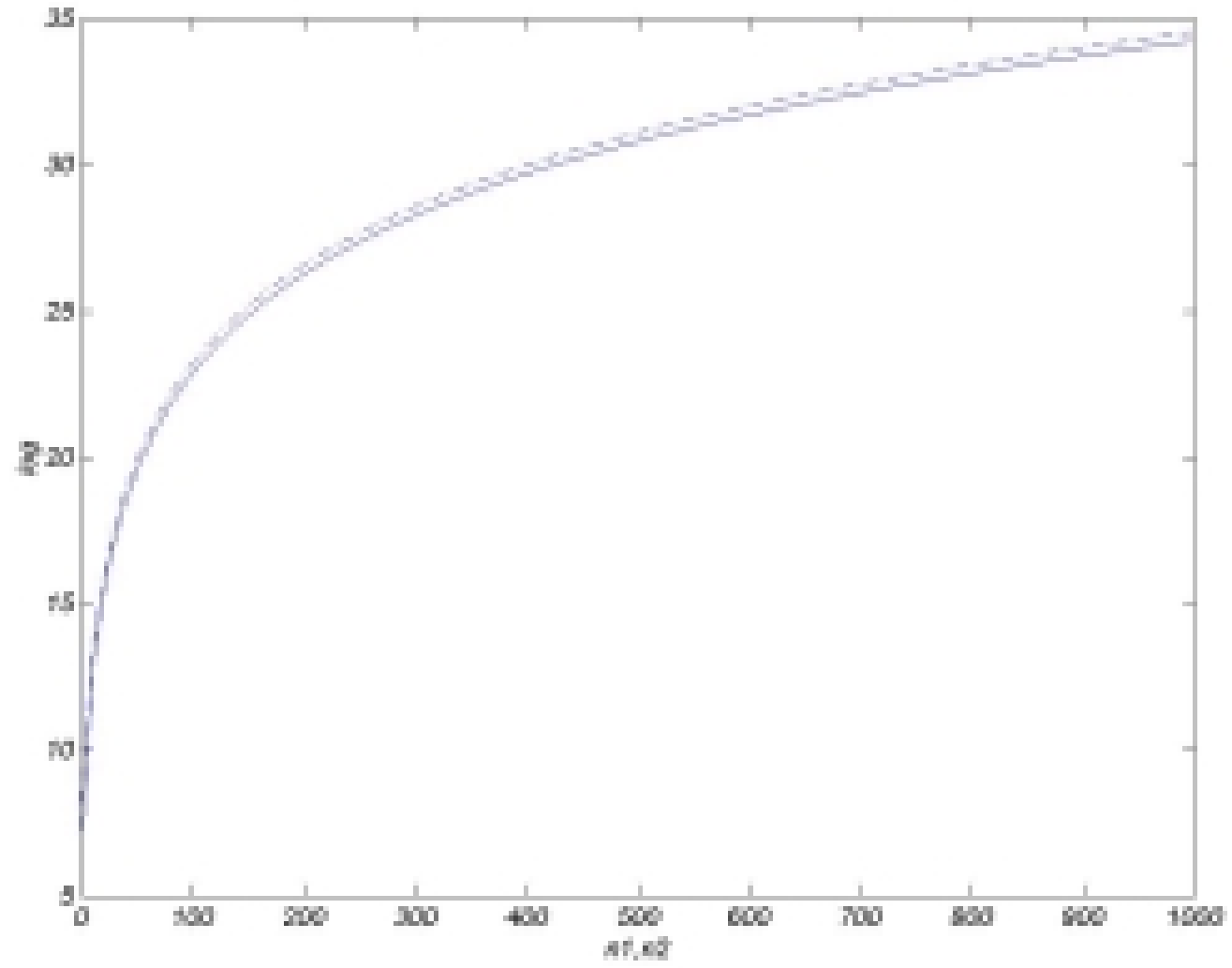

*Graph 1*

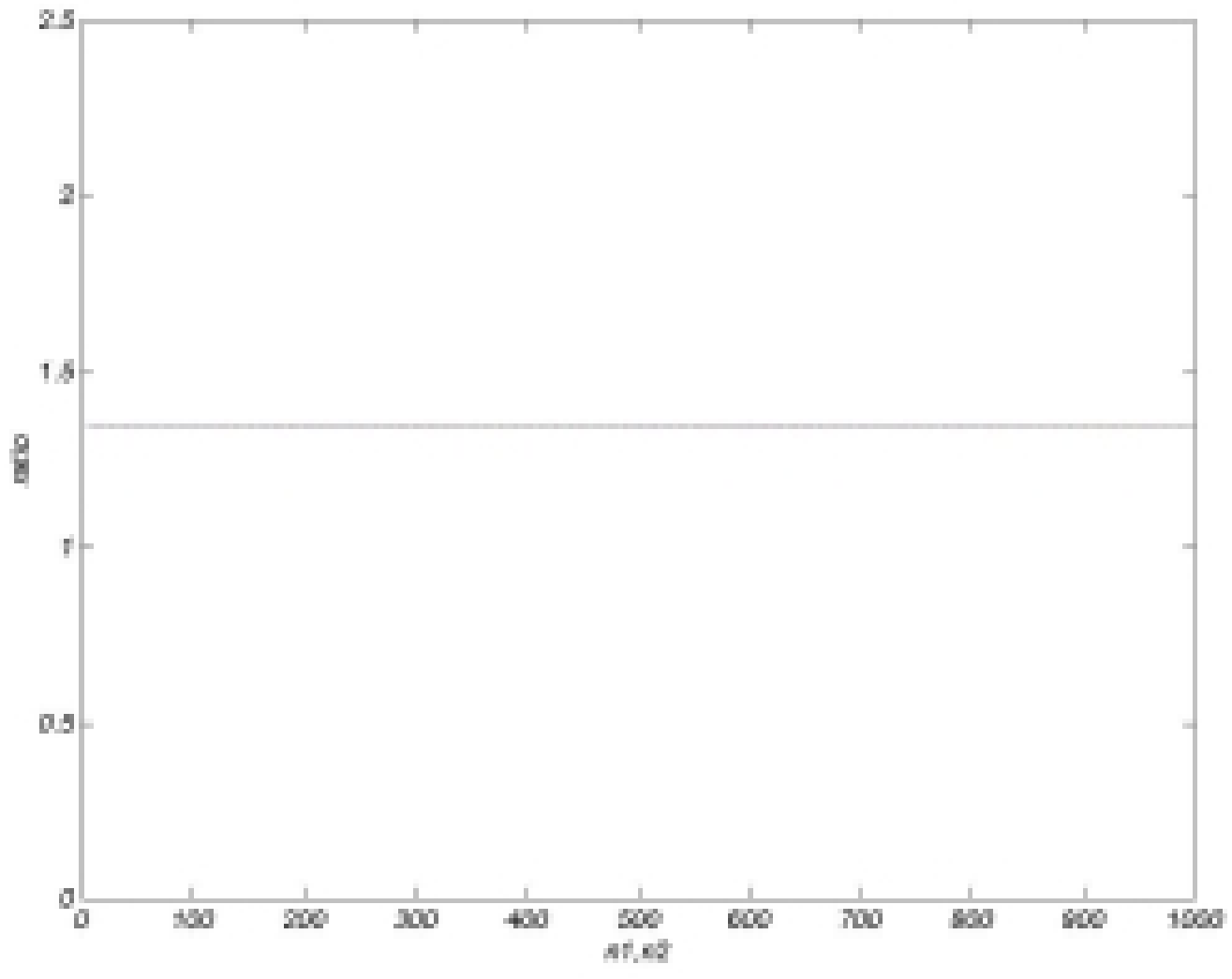

*Graph 2*

**Conclusion.** The previous presented results, explain clearly the structure of the Generalized Reflection coefficients, and their relevant orthogonal polynomials in the Hermetian Toeplitz Block Toeplitz Case, while the reflection coefficients admit the same Block Toeplitz recurrence, and the polynomials $\mathrm{p}_{k,l}$ $\mathrm{q}_{k,l}$, show a block Toeplitz recurrence with an added shift between Blocks polynomials counterparts; both of reflection coefficients, and polynomials $a_{k,l}, a'_{k,l}, \mathrm{p}_{k,l}, \mathrm{q}_{k,l}$ with indices $k,l$, shows an exchange structure with values of indices $k \otimes l$, symmetric to the antidiagonal of each block. The later fact did establish a faster, levinson-like inversion algorithm.

## References.